\documentclass[a4paper,12pt]{article}
\usepackage{amsmath,amsfonts,calrsfs,amssymb,color}
\newtheorem{Theorem}{Theorem}[section]
\newtheorem{Definition}[Theorem]{Definition}
\newtheorem{Proposition}[Theorem]{Proposition}
\newtheorem{Lemma}[Theorem]{Lemma}

\newtheorem{Remark}[Theorem]{Remark}

\newcommand{\one}{1\!\!\!\;\mathrm{l}}

\makeatletter
\newif\ifmsbmloaded@

\@addtoreset{equation}{section}

\makeatother

\def\R{\mathbb R}
\def\N{\mathbb N}

\def\E{\mathbb E}
\def\P{\mathbb P}
\def\ds{\displaystyle}

\title{\bf On the martingale probem associated to the $2D$ and $3D$ Stochastic Navier-Stokes equations}
\author{Giuseppe Da Prato
\\Scuola Normale Superiore di Pisa,\\  Piazza dei Cavalieri
7, 56126 Pisa, Italy
\and
Arnaud Debussche\\
Ecole Normale Sup\'erieure de Cachan,\\
 antenne de
Bretagne,\\
 Campus de Ker Lann, 35170 Bruz, France
}
\date{}
\begin{document}
\maketitle

\begin{abstract} 
We consider the martingale problem associated to the Navier-Stokes in dimension $2$ or $3$.  Existence is well known and it has been recently shown that markovian transition semi group associated 
to these equations can be constructed. We study the Kolmogorov operator associated to these 
equations. It can be defined formally as a differential  operator on an infinite dimensional Hilbert space. 
It can be also defined in an abstract way as the infinitesimal generator of the transition semi group. We 
explicit cores for these abstract operators and identify them with the concrete differential operators on 
these cores. In dimension $2$, the core is explicit and we can use a classical argument to prove uniqueness for the martingale problem. In dimension $3$, we are only able to exhibit a core which is defined abstractly and does not allow to prove uniqueness for the martingale problem. Instead, we exhibit a core for a modified Kolmogorov operator which enables us to prove uniqueness for the martingale problem up to the time the solutions are regular. 
\end{abstract}
\bigskip

\noindent {\bf 2000 Mathematics Subject Classification AMS}: %76D05, 60H15,37A25   

\noindent {\bf Key words}: Stochastic Navier--Stokes, Kolmogorov equations, martingale problems, 
weak uniqueness.
\bigskip

\section{Introduction}

We consider the stochastic Navier--Stokes on a bounded domain
${\mathcal O}$ of $\R^d$, $d=2$ or $3$, with 
Dirichlet  boundary conditions:
the unknowns are the velocity $X(t,\xi)$ and the pressure
 $p(t,\xi)$ defined for $t>0$ and $\xi\in \overline{\mathcal O}$:
\begin{equation}
\label{e1.0} \left\{\begin{array}{l} dX(t,\xi)=[\Delta
X(t,\xi)-(X(t,\xi)\cdot\nabla)X(t,\xi)]dt-\nabla
p(t,\xi)dt+f(\xi)dt+\sqrt{Q}\;dW,
\\
\mbox{\rm div}\;X(t,\xi)=0,
\end{array}\right.
\end{equation}
with Dirichlet boundary conditions
$$
X(t,\xi)=0,\quad t>0,\;\xi\in \partial \mathcal O,
$$
and supplemented with the initial condition
$$
X(0,\xi)=x(\xi),\;\xi\in  \mathcal O.
$$
We have taken the viscosity equal to $1$ since it plays no particular role in this work.  The 
 understanding of the  stochastic Navier-Stokes equations have progressed considerably 
 recently. In dimension two, impressive progresses have been obtained
and difficult ergodic properties have been proved (see\cite{kupiainen}, \cite{mattingly}, 
 \cite{flandoli-maslowski}, \cite{HM},
\cite{HM2}, \cite{kuksinBook}, \cite{kuksin}, \cite{kuksin2}, \cite{mattingly2},
\cite{mattingly3}). In dimension three, the theory is not so advanced. Uniqueness is still an
open problem. However, Markov solutions have been constructed and 
ergodic properties have been proved recently (see \cite{chueschov-kuksin}, \cite{chueschov-kuksin2}, \cite{DPD-NS3D}, \cite{DO06}, 
\cite{flandoli-cetraro}, \cite{romito-flandoli-selection}, \cite{Odasso}, \cite{romito}, \cite{shirikyan3D}).

In this article, our aim is to try to improve the understanding of the martingale problems associated
to these equations. Let us first set some notations. Let
$$
H=\{x\in (L^2(\mathcal O))^d:\;\mbox{\rm div}\; x=0\;\mbox{\rm
in}\;\mathcal O,\;x\cdot n=0\;
\mbox{\rm on}\;\partial \mathcal O\},
$$
where $n$ is  the outward normal to $\partial \mathcal O,$
and $V= (H^1_0(\mathcal O))^d\cap H.$
The norm and  inner product in $H$   will be denoted by $|\cdot|$   and
 $(\cdot,\cdot)$ respectively. Moreover $W$ is a cylindrical Wiener process  on $H$ and the 
 covariance of the noise $Q$ is trace class and non degenerate (see \eqref{e1.2} and 
 \eqref{e1.3} below for more precise assumptions).

  \vspace{2mm}
 \noindent We also denote by $A$ the Stokes operator in $H$:
    $$A= P\Delta,\quad D(A)=(H^2(\mathcal O))^d\cap V,$$
where   $P$ is the orthogonal projection of $(L^2(\mathcal O))^3$ onto $H$
and by $b$  the operator
$$
b(x,y) =-P((x\cdot\nabla)y),\quad
b(x)=b(x,x),\quad x,y\in V.
$$

With these notations we rewrite the equations as
\begin{equation}
\label{e1.1}
\left\{\begin{array}{l}
dX=(AX+b(X))dt+\sqrt Q\;dW,\\
\\
X(0)=x.
\end{array}\right. 
\end{equation}
We assume that
\begin{equation}
\label{e1.2}
\quad\mbox{\rm Tr}\;(-A)^{1+g}Q<\infty,\quad \mbox{\rm for some}\;g>0 
\end{equation}
 and
\begin{equation}
\label{e1.3}
|Q^{-1/2}x|\le c|(-A)^rx|,\quad \mbox{\rm for some}\;r\in (1,3/2).
\end{equation} 
In dimension $d=3$, it is well known that there exists a solution to the martingale problem but 
weak or strong uniqueness is an open problem (see \cite{flandoli-cetraro} for a survey).  
However, it has been proved in \cite{DPD-NS3D}, \cite{DO06} (see also \cite{romito-flandoli-selection}) that the above assumptions allow to construct a transition semigroup $(P_t)_{t\ge 0}$ associated to a Markov family of solutions
$$
((X(t,x)_{t\ge 0},\Omega_x,\mathcal F_x,\P_x)
$$
 for $x\in D(A)$. Moreover for sufficiently regular $\varphi$ defined on $D(A)$, $P_t\varphi$ is a solution of
 the Kolmogorov equation associated to \eqref{e1.1}
 \begin{equation}
\label{e1.4}
\left\{\begin{array}{l}
\ds \frac{du}{dt}=Lu, \; t>0, \; x\in D(A),\\
\\
u(0,x )=\varphi(x),\; x\in D(A),
\end{array}\right. 
\end{equation}
where the Kolmogorov operator $L$ is defined 
by 
$$
L\varphi(x)=\frac12 {\rm Tr}\left\{ Q D^2 \varphi(x) \right\} + (Ax+b(x), D\varphi(x))
$$
for sufficiently smooth functions $\varphi$ on $D(A)$.

In all the article, we choose one Markov family $((X(t,x)_{t\ge 0},\Omega_x,\mathcal F_x,\P_x)$ as the one 
constructed in \cite{DO06}.

The fundamental idea in \cite{DPD-NS3D} is to introduce a modified semigroup $(S_t)_{t\ge 0}$ defined by
\begin{equation}
\label{e1.4bis}
S_t\varphi(x)= \E(e^{-K\int_0^t |AX(s,x)|^2 ds } \varphi(X(t,x))) .
\end{equation}
It can be seen that for $K$ large enough,
this semigroup has very nice smoothing properties and various estimates 
can be proved. Note that, thanks to Feynman-Kac formula,  this semigroup is formally
associated to the following equation
 \begin{equation}
\label{e1.5}
\left\{\begin{array}{l}
\ds \frac{dv}{dt}=Nv, \; t>0, \; x\in D(A),\\
\\
v(0,x )=\varphi(x),\; x\in D(A),
\end{array}\right. 
\end{equation}
where $N$ is defined
$$
N\varphi(x) = \frac12 {\rm Tr}\left\{ Q D^2 \varphi(x) \right\} + (Ax+b(x), D\varphi(x))-K|Ax|^2\varphi(x),
$$
for sufficiently smooth functions $\varphi$ on $D(A)$.

In \cite{DPD-NS3D}, \cite{DO06}, this semigroup is defined only on the Galerkin approximations of \eqref{e1.1}.
Let $P_m$ denote the projector associated to the first $m$ eigenvalues of $A$. We consider the 
following equation in $P_mH$
 \begin{equation}
\label{e1.6}
\left\{\begin{array}{l}
dX_m=(AX_m+b_m(X_m))dt+\sqrt Q_m\;dW\\
\\
X_m(0)=P_mx,
\end{array}\right. 
\end{equation}
where $b_m(x) = P_m b(P_m x)$, $Q_m =P_m QP_m$. This defines, with obvious notations, 
$(P_t^m)_{t\ge 0}$ and $(S_t^m)_{t\ge 0}$. The following formula holds by a standard argument:
$$
P_t^m\varphi = S_t^m \varphi + K\int_0^t S_{t-s}^m 
\left( |A\cdot|^2 P^m_s\varphi ds\right),\; \varphi\in C_b(P_mH).
$$
Various estimates are proved on $(S_t^m)_{t\ge 0}$ and transferred to $(P_t^m)_{t\ge 0}$ thanks to this identity. A compactness argument allows to construct $(P_t)_{t\ge 0}$. Moreover, 
a subsequence $m_k$ can be constructed such that for any $x\in D(A)$, $(X_{m_k}(\cdot,x))_{t\ge 0}$ 
converges in law to $(X(\cdot,x))_{t\ge 0}$.
  
Note also that similar arguments as in \cite{DPD-NS3D} 
may be used to prove that for smooth $\varphi$, $(S_t\varphi)_{t\ge 0}$ is a strict solution to 
\eqref{e1.5}. 

 In dimension $2$ this result also holds with exactly the same proofs since all arguments for $d=3$ are still valid.
 Note that  it is well known that for $d=2$ conditions \eqref{e1.2}-\eqref{e1.3} imply that, for $x\in H$, there exists a unique strong solution to \eqref{e1.1} and the proof of the above facts can be simplified. 
  
 In the following, we give some properties of the generator of $(P_t)_{t\ge 0}$ and  $(S_t)_{t\ge 0}$.
 For $d=2$, we explicit a core, identify the abstract generator with the differential 
 operator $L$ on this core and prove existence and uniqueness for the corresponding martingale problem. 
 (See \cite{roeckner-sobol}Ê  for a similar result). Again, this follows from strong uniqueness but we think that it is interesting to have a 
 direct proof of this fact. Moreover, it can be very useful to have a better knowledge of the Kolmogorov 
 generator and we think that this work is a contribution in this direction. In dimension $3$, 
 we are not able to prove this.  
 We explain the difficulties encountered. We hope that this article will help the reader to get 
 a better insight into the problem of weak uniqueness for the three dimensional Navier-Stokes equations.
Nonetheless, we explicit a core for the generator of the transformed semigroup $(S_t)_{t\ge 0}$,
identify it with the differential operator $N$ on this core
and prove uniqueness for the stopped martingale problem. In other words, we prove weak uniqueness
up to the time solutions are smooth. Again, this could be proved directly thanks to local strong 
uniqueness.

\section{The generators}

The space of continuous functions on $D(A)$ is denoted by $C_b(D(A))$. Its norm is
denoted by $\|\cdot\|_0$. For $k\in \N$, $C^k(D(A))$ is the space of $C^k$ functions on $D(A)$.
We need several other function spaces on $D(A)$.

Let us introduce the set $\mathcal E_1\subset  C_b(D(A))$  of $C^3$ functions on $D(A)$ such that
there exists a constant $c$ satistying
\begin{itemize}
\item $|(-A)^{-1} Df(x)|_{H} \le c (|Ax|^2+1)$
\item $|(-A)^{-1}D^2f(x)(-A)^{-1}|_{\mathcal L(H)} \le c (|Ax|^4+1)$
\item $ |(-A)^{-1/2}D^2f(x)(-A)^{-1/2}|_{\mathcal L(H)} \le c (|Ax|^6+1)$
\item $ \|D^3f(x)\left((-A)^{-1}\cdot,(-A)^{-1}\cdot,(-A)^{-1}\cdot\right)\|  \le c (|Ax|^6+1) $
\item $ \|D^3f(x)\left((-A)^{-\gamma}\cdot,(-A)^{-\gamma}\cdot,(-A)^{-\gamma}\cdot\right)\|  \le c (|Ax|^8+1) $
\item $|Df(x)|_{H} \le c (|Ax|^4+1)$
\end{itemize}
where $\gamma\in (1/2,1]$ and
$$
\begin{array}{l}
\ds \mathcal E_2= \left\{f\in  C_b(D(A)),\; \sup_{x,y\in D(A)}\frac{|f(x)-f(y)|}
{|A(x-y)|(1+|Ax|^2+|Ay|^2}<+\infty\right\}.
\end{array}
$$
Note that we identify the gradient and the differential of a real valued function. Also, the second differential 
is identified with a function with values in $\mathcal L(H)$. The third differential is a trilinear operator
on $D(A)$ and the norm $\|\cdot\|$ above is the norm of such operators.

Slightly improving the arguments in \cite{DPD-NS3D}, it can be  proved{\footnote{In fact, only 
Lemma 5.3 has to be improved. In this Lemma, the term $L_1$ can in fact 
be estimated in a single step by using Proposition 3.5 of \cite{DO06} instead of Proposition 
5.1 of  \cite{DPD-NS3D}. }} that $P_t$ 
maps $\mathcal E_i$ into itself
and that there exists a constant $c>0$ such that
\begin{equation}
\label{e4}
\|P_t f\|_{\mathcal E_i}\le c\|f\|_{\mathcal E_i}.
\end{equation}
Moreover, for $f\in \mathcal E_1$, $P_t f$ is a strict solution of \eqref{e1.4} in the sense that it is satisfied for any $x\in D(A)$ and $t\ge 0$. Again, the result of  \cite{DPD-NS3D} has to be slightly improved to get this
result. In fact, using an interpolation argument, Proposition 5.9 and the various other estimates  in 
\cite{DPD-NS3D}, it is easy to deduce that, for any $x\in D(A)$,  $LP_tf(x)$ is continuous on 
$[0,T]$. 

For $f\in \mathcal E_2$, $P_tf$ is still a solution of \eqref{e1.4}  but in the mild sense. We define the 
Ornstein-Uhlenbeck semigroup associated to the linear equation
$$
R_t \varphi (x) = \varphi( e^{tA}x + \int_0^t e^{A(t-s)}\sqrt{Q} dW(s), \; t\ge 0, \; \varphi \in C_b(D(A)).
$$
Then it is shown in \cite{DPD-NS3D} that
 \begin{equation}
\label{e2.1bis}
P_tf(x)= R_tf(x) +\int_0^t R_{t-s}(b, DP_sf)ds,\; t\ge 0,\;  f\in \mathcal E_2.
\end{equation}
For any $\lambda>0$ we set
$$
F_\lambda f=\int_0^\infty e^{-\lambda t}P_t f dt,\quad f\in C_b(D(A)).
$$
Then since
$\|P_tf\|_0\le \|f\|_0$, we have
$$
\|F_\lambda f\|_0\le \frac1\lambda\|f\|_0.
$$
Moreover, since $P_t$ is Feller, we have by dominated convergence
$$
F_\lambda f\in C_b(D(A)).
$$
It can be easily deduced that 
$$
F_\lambda f-F_\mu f=(\mu-\lambda)F_\lambda F_\mu f,\quad  \mu,\lambda>0,
$$
and
$$
\lim_{\lambda\to \infty}\lambda F_\lambda f(x)=
\lim_{\lambda\to \infty}\int_0^\infty e^{-\tau}P_{\tau/\lambda}(x)d\tau=f(x),\quad x\in D(A).
$$
It follows classically (see for instance \cite{priola}) that there exists a unique maximal dissipative operator $\bar L$ on $C_b(D(A))$ with domain $D(\bar L)$ such that
$$
F_\lambda f=(\lambda-\bar L)^{-1} f.
$$
We recall the following well known characterization of $D(\bar L)$:
$f\in D(\bar L)$ if and only if
\begin{enumerate}
\item[(i)] $f\in C_b(D(A)),$

\item[(ii)] $\ds \frac1t\;\|P_t f-f\|_0$ is bounded for $t\in[0,1]$,

\item[(iii)] $\ds \frac1t\;(P_t f(x)-f(x))$ has a limit for any $x\in D(A)$.
\end{enumerate}\bigskip
Moreover, we have in this case
$$
\bar Lf(x)=\lim_{t\to 0}\frac1t\;(P_tf(x)-f(x)).
$$
Recall also that
$$
(\lambda-\bar L)^{-1} f=\int_0^\infty e^{-\lambda t}P_t f dt,\quad f\in 
C_b(D(A)).
$$
By \eqref{e4} we deduce that
\begin{equation}
\label{e2.2}
\|(\lambda-\bar L)^{-1} f\|_{\mathcal E_i}\le \frac{c}\lambda\|f\|_{\mathcal E_i}.
\end{equation}

Similarly, we may define, for $k\ge 0$, $\mathcal E_3^k$ as the space $C^3$ functions on $D(A)$ such that
there exists a constant $c$ satistying
\begin{itemize}
\item $|(-A)^{-1} Df(x)|_{H} \le c (|Ax|^k+1)$
\item $|(-A)^{-1}D^2f(x)(-A)^{-1}|_{\mathcal L(H)} \le c (|Ax|^k+1)$
\item $ |(-A)^{-1/2}D^2f(x)(-A)^{-1/2}|_{\mathcal L(H)} \le c (|Ax|^k+1)$
\item $ \|D^3f(x)\left((-A)^{-1}\cdot,(-A)^{-1}\cdot,(-A)^{-1}\cdot\right)\|  \le c (|Ax|^k+1) $
\item $ \|D^3f(x)\left((-A)^{-\gamma}\cdot,(-A)^{-\gamma}\cdot,(-A)^{-\gamma}\cdot\right)\|  \le c (|Ax|^k+1) $
\item $|Df(x)|_{H} \le c (|Ax|^k+1)$
\end{itemize}
where $\gamma\in (1/2,1]$. By the various estimates given in \cite{DPD-NS3D}, it is easy to check that,
provided $K$ is chosen large enough, 
$S_t$ maps $\mathcal E_3^k$ into itself and 
there exists a constant $c>0$ such that
\begin{equation}
\label{e4bis}
\| S_t f\|_{\mathcal E_3^k}\le c\|f\|_{\mathcal E_3^k}.
\end{equation}
Moreover, for $f\in \mathcal E_3^k$, $S_t f$ is a strict solution of \eqref{e1.5} in the sense that it is satisfied for any $x\in D(A)$ and $t\ge 0$. 

For any $\lambda>0$ we set
$$
{\widetilde F}_\lambda f=\int_0^\infty e^{-\lambda t}S_t f dt,\quad f\in C_b(D(A)).
$$
and prove that there exists a unique maximal dissipative operator $\bar N$ on $C_b(D(A))$ with domain $D(\bar N)$ such that
$$
{\widetilde F}_\lambda f=(\lambda-\bar N)^{-1} f,
$$
and
$f\in D(\bar N)$ if and only if
\begin{enumerate}
\item[(i)] $f\in C_b(D(A)),$

\item[(ii)] $\ds \frac1t\;\|S_t f-f\|_0$ is bounded for $t\in[0,1]$,

\item[(iii)] $\ds \frac1t\;(S_t f(x)-f(x))$ has a limit for any $x\in D(A)$.
\end{enumerate}\bigskip

Finally, by \eqref{e4bis}, we see that 
\begin{equation}
\label{e2.2bis}
\|(\lambda-\bar N)^{-1} f\|_{\mathcal E_3^k}\le \frac{c}\lambda\|f\|_{\mathcal E_3^k}.
\end{equation}

\section{Construction of cores  and identification of
the generators }

In this section, we analyse the generators defined in the preceeding section. We start with the 
following definition.
\begin{Definition}  Let $K$ be an operator with domain $D(K)$. A set $\mathcal D\subset D(K)$ is a 
$\pi$-core for $K$ if for 
any $\varphi \in D(K)$, there exists a sequence $(\varphi_n)_{n\in\N}$ in
$\mathcal D$ which $\pi$-converges{\footnote{Recall that the $\pi$-convergence - also called 
b.p. convergence - is defined by : $(f_n)_{n\in\N}$ $\pi$-converges to $f$ iff $f_n(x)\to f(x)$ for
any $x\in D(A)$ and $\sup_{n\in\N}\|f_n\|_0<\infty$.}}
to $\varphi$ and such that $(K\varphi_n)_{n\in\N}$ $\pi$-converges 
to $K \varphi$.
\end{Definition}

Let us set $\mathcal G_1=(\lambda-\bar L)^{-1}\mathcal E_1$ for some $\lambda>0$. Clearly for any $\varphi\in \mathcal G_1$ we have
$\varphi\in D(\bar L)$ and by \eqref{e2.2}, $\varphi\in \mathcal E_1.$ Moreover,
$$
P_t\varphi(x)-\varphi(x)=\int_0^t LP_s\varphi(x)ds,
$$
since $(P_t\varphi)_{t\ge 0}$ is a strict solution of the Kolmogorov equation. By \eqref{e4} and the definition
of $\mathcal E_1$, for any $x\in D(A)$ we have
$$
|LP_s\varphi(x)|\le c(1+|Ax|^6)\|P_s\varphi\|_{\mathcal E_1}
\le c(1+|Ax|^6)\|\varphi\|_{\mathcal E_1}.
$$
Moreover, since
$$
t\to LP_t\varphi(x)
$$
is continuous, we have
$$
\lim_{t\to 0}\frac1t\;(P_t\varphi(x)-\varphi(x))=L\varphi(x).
$$
We deduce that
$$
\bar L\varphi(x)=L\varphi(x),\quad x\in D(A).
$$
Since $\mathcal E_1$ is $\pi$-dense in $C_b(D(A))$, we deduce that
$\mathcal G_1$ is a $\pi$-core for $\bar L$.

Also $\mathcal E_1\subset \mathcal E_2$ so that
$$
\mathcal G_1\subset \mathcal G_2=(\lambda-\bar L)^{-1}\mathcal E_2
$$
and $\mathcal G_2$ is also a $\pi$-core for $\bar L$.

These results hold both in dimension $2$ or $3$. The problem is that these cores are abstract
and strongly  depend 
on the semigroup $(P_t)_{t\ge 0}$. In dimension $3$, this is a real problem since we do not 
know if the transition semigroup is unique. If we were able to construct a core in terms of the differential
operator $L$, this would certainly imply uniqueness of this transition semigroup.

In dimension $2$, we are able to construct such a core. Of course, in this case, uniqueness
is well known. However, we think that it is important to have explicit cores. This gives many informations
on the transition semigroup $(P_t)_{t\ge 0}$. 

\begin{Theorem}
\label{t3.1}

Let us set
$$
\mathcal H=\{f\in \mathcal E_1:\;Lf\in \mathcal E_1\}
$$
then, in dimension $d=2$,  $\mathcal H\subset D(\bar L)$ and it is a $\pi$-core for $\bar L$. Moreover, for any 
$f \in \mathcal H$, we have 
$$
\bar L f = L f.
$$
\end{Theorem}

The crucial point is to prove the following result.
\begin{Proposition}
\label{p3.2}
Let $d=2$. For any  $f\in \mathcal H$
we have
$$
P_{t_1}f-P_{t_2}f=\int_{t_1}^{t_2}P_sLfds,\quad0\le t_1\le t_2.
$$
\end{Proposition}
{\bf Proof}. Let $f\in \mathcal H$. By It\^o
formula applied to the Galerkin equation \eqref{e1.6}, we have for $\epsilon >0$
$$
\begin{array}{l}
\ds d\left( e^{-\epsilon\int_0^t|(-A)^{1/2}X_m(s,x)|^6 ds}f(X_m(t,x))  \right) 
\\
\\
\ds=\left(-\epsilon|(-A)^{1/2}X_m(t,x)|^6 f(X_m(t,x))+L_mf_m(X(t,x))\right)
e^{-\epsilon\int_0^t|(-A)^{1/2}X_m(s,x)|^6 ds} dt\\
\\
\ds +e^{-\epsilon\int_0^t|(-A)^{1/2}X_m(s,x)|^6 ds}(Df_m(X_m(t,x)),\sqrt Q_m\;dW)
\end{array}
$$
and
\begin{equation}
\label{e5}
\begin{array}{l}
\ds \E \left( e^{-\epsilon\int_0^{t_2}|(-A)^{1/2}X_m(s,x)|^6 ds}f(X_m(t_2,x))  \right)\\
\\
\ds-
\left( e^{-\epsilon\int_0^{t_1}|(-A)^{1/2}X_m(s,x)|^6 ds}f(X_m(t_1,x))  \right)
\\
\\
\ds =\E\Big(\int_{t_1}^{t_2}\big(-\epsilon|(-A)^{1/2}X_m(s,x)|^6 f(X_m(s,x))\\
\\
\ds
+L_mf(X_m(s,x))\big)e^{-\epsilon\int_0^{s}|(-A)^{1/2}X_m(\sigma,x)|^6 d\sigma}ds \Big). 
\end{array} 
\end{equation}
We have denoted by $L_m$ the Kolmogorov operator associated to \eqref{e1.6}. 
Since $f\in \mathcal H$, we have 
$$
|Lf_m (x)|\le c(1+|Ax|^6).
$$
By Proposition \ref{p3.1} and Lemma \ref{l2bis}, the right hand side of
\eqref{e5} is uniformly integrable on $\Omega\times [t_1,t_2]$ with respect to $m$.
Thus, we can take the limit $m\to \infty$ in \eqref{e5} and obtain
\begin{equation}
\label{e5bis}
\begin{array}{l}
\ds \E_x \left( e^{-\epsilon\int_0^{t_2}|(-A)^{1/2}X(s,x)|^6 ds}f(X(t_2,x))  \right)\\
\\
\ds-
\E_x\left( e^{-\epsilon\int_0^{t_1}|(-A)^{1/2}X(s,x)|^6 ds}f(X(t_1,x))  \right)
\\
\\
\ds =\E_x\Big(\int_{t_1}^{t_2}\big(-\epsilon|(-A)^{1/2}X(s,x)|^6 f(X(s,x))\\
\\
\ds
\hspace{5mm}+Lf(X(s,x))\big)e^{-\epsilon\int_0^{s}|(-A)^{1/2}X(\sigma,x)|^6 d\sigma}ds \Big). 
\end{array} 
\end{equation}

\bigskip It is easy to prove by dominated convergence that
$$
\begin{array}{l}
\ds\E_x\left( e^{-\epsilon\int_0^{t_i}|X(s,x)|^6_1ds}f(X(t_i,x))  \right)\to P_{t_i}f(x),
\\
\\
\ds\E_x\int_{t_1}^{t_2}Lf(X(s,x))e^{-\epsilon\int_0^{s}|X(\sigma,x)|^6_1d\sigma}ds\to \E_x\int_{t_1}^{t_2}P_sLf(x),
\end{array}
$$
when $\epsilon\to 0$. Indeed by Lemma \ref{l2bis} below, we have
$$
\int_0^{t_i}|X(s,x)|^6_1ds<\infty\quad \P\mbox{\rm -a.s.}.
$$
Moreover
$$
\begin{array}{l}
\ds\left|\E_x \int_{t_1}^{t_2}\left(\epsilon|X(s,x)|^6_1f(X(s,x))
e^{-\epsilon\int_0^s|X(\sigma,x)|^6_1d\sigma} ds  \right)    \right|
\\
\\
\le \|f\|_0\E_x \left(
e^{-\epsilon\int_0^{t_1}|X(\sigma,x)|^6_1d\sigma}-
e^{-\epsilon\int_0^{t_2}|X(\sigma,x)|^6_1d\sigma}    \right)\to 0,
\end{array}
$$
as $\epsilon\to 0$. The result follows. $\Box$

It is now easy to conclude the proof of Theorem \ref{t3.1}. Indeed, by Proposition \ref{p3.2}, for
$f\in \mathcal H$ we have, since $\|P_s Lf\|_0\le \|Lf\|_0$, 
$$
\|P_t f -f\|_0 \le t \| Lf\|_0.
$$
Moreover, since $s\mapsto P_sLf(x)$, is continuous for any $x\in D(A)$
$$
\frac1t \left(P_t f(x) -f(x)\right) \to Lf(x),\mbox{ as } t\to 0.
$$
It follows that $f\in D(\bar L)$ and $\bar L f = Lf$. 
Finally
$$
\mathcal G_1\subset  \mathcal H
$$
and since $\mathcal G_1$ is a $\pi$-core we deduce that $\mathcal H$ is also a $\pi$-core. 

\begin{Remark}
\em We do not use that $P_tf$ is a strict solution of the Kolmogorov equation to prove 
that $\mathcal H\subset D(\bar L)$ and $\bar L f=Lf$.
But we do not know if there is a direct proof of the fact that
$\mathcal H $ is a $\pi$-core. We have used that $\mathcal G_1\subset \mathcal H$ and that $\mathcal G_1 $ is a $\pi$-core.
The proof of $\mathcal G_1 \subset\mathcal H$ requires \eqref{e4}
which is almost as strong as the construction of a strict solution. $\Box$
\end{Remark}
\begin{Remark}
\em For $d=3$, using Lemma 3.1 in \cite{DPD-NS3D},  it is easy to prove a formula similar to \eqref{e5bis} with
$|(-A)^{1/2} X(s,x)|^6$ replaced by $|AX(s,x)|^4$ in the exponential terms. The problem is that
$$
\lim_{\epsilon\to 0}e^{-\epsilon\int_0^{t}|AX(s,x)|^4 ds}=\one_{[0,\tau^*(x))}(t),
$$
where $\tau^*(x)$ is the life time of the solution in $D(A)$. Thus we are not able to prove Proposition 
\ref{p3.2} in this case. 
\end{Remark}
We have the following result on the operator $\bar N$.
\begin{Theorem}
\label{t3.2}
Let $d=2$ or $3$ and $k\in \N$, define
$$
\widetilde{\mathcal H_k}= \{f\in \mathcal E_3^k\; : Nf\in \mathcal E_3^k\}.
$$
Then $\widetilde{\mathcal H_k}\subset D(\bar N)$ and it is $\pi$-core for $\bar N$. Moreover, for any 
$f\in \widetilde H_k$ we have 
$$
\bar N f = Nf.
$$
\end{Theorem}

The proof follows the same line as above. Indeed, it is easy to use similar arguments as in \cite{DPD-NS3D}
and prove that for $f\in  \mathcal E_3^k$, $(S_tf )_{t\ge 0}$ is a strict solution to \eqref{e1.5}. Arguing as
above, we deduce that $(\lambda-\bar N)^{-1} \mathcal E_3^k$ is a $\pi$-core for $\bar N$. 

Moreover, applying It\^o formula to the Galerkin approximations and letting $m\to \infty$ along the 
subsequence $m_k$ - thanks to 
Lemma 3.1 of \cite{DPD-NS3D} to get uniform integrability - we prove, for $f\in \widetilde{\mathcal H_k}$,
$$
\begin{array}{l}
\ds \E_x(e^{-K\int_{0}^{t_2} |AX(s,x)|^2ds} f(X(t_2,x) )-\E_x(e^{-K\int_{0}^{t_1} |AX(s,x)|^2ds} f(X(t_1,x) )\\
\\
\ds = \E_x \left( \int_{t_1}^{t_2} \left( -K|AX(s,x)|^2 f(X(s,x)) + L f(X(s,x)) \right) e^{-K\int_{0}^{s} |AX(\sigma,x)|^2
d\sigma} ds\right).
\end{array}
$$
We rewrite this as 
$$
S_{t_2}f(x)-S_{t_1}f(x)=\int_{t_1}^{t_2} S_s Nf (x),
$$
and deduce as above that $f\in D(\bar N)$ and $\bar N f =Nf$. Finally, since
$(\lambda-\bar N)^{-1} \mathcal E_3^k \subset \widetilde H_k$, we know that $\widetilde{\mathcal H_k}$ is also 
a $\pi$-core.

\section{Uniqueness for the martingale problem}

Let us study the following martingale problem.
\begin{Definition}
\label{d1}
We say that a probability measure $\P_x$  on $C([0,T];D((-A)^{-\epsilon}))$, $\epsilon>0$ is a solution of the martingale problem
associated to \eqref{e1.1}  if
$$
\P_x(\eta(t)\in D(A))=1,\quad t\ge 0, \; \P_x(\eta(0)=x)=1
$$
and for any $f\in \mathcal H$
$$
f(\eta(t))-\int_0^tLf(\eta(s))ds,
$$
is a martingale with respect to the natural  filtration.
\end{Definition}

\begin{Remark}
\em  In general, it is proved the existence of a solution to a different   martingale problem where $f$ is required to be in a smaller class. In particular, it is required that $f\in C_b(D((-A)^{-\epsilon}))$ for some  $\epsilon>0$. However, in all concrete construction of solutions,  it can be shown that a solution of our 
martingale problem is in fact obtained.  $\Box$
\end{Remark}
\begin{Theorem}
\label{t4.1}
Let $d=2$, then for any $x\in D(A)$, there exists a unique solution to the martingale problem.
\end{Theorem}
{\bf Proof.} 
By a similar proof as for  Proposition \ref{p3.2}, we know that there exists a solution to the martingale problem

Uniqueness follows from a classical argument. 
Let $f\in \mathcal E_1$ and, for $\lambda>0$ set $\varphi=(\lambda-\bar L)^{-1}f.$ Then $\varphi\in \mathcal G_1\subset \mathcal H$ and
$$
\varphi(\eta(t))-\varphi(x)-\int_0^t L\varphi(\eta(s))ds
$$
is a martingale. Thus, for any solution $\tilde P_x$ of the martingale problem,
$$
\tilde \E_x\left(\varphi(\eta(t))-\varphi(x)-\int_0^tL\varphi(\eta(s))ds\right)=\varphi(x).
$$
We multiply by $\lambda e^{-\lambda t}$, integrate over $[0,\infty)$ and obtain, since 
$\bar L\varphi =L\varphi$,
$$
\tilde\E_x\int_0^\infty e^{-\lambda t}f(\eta(t))dt=\varphi(x)=(\lambda-\bar L)^{-1}f(x)=\int_0^\infty e^{-\lambda t}P_t f(x)dt.
$$
By   inversion of Laplace transform we deduce
$$
\tilde\E_x(f(\eta(t))=P_tf(x).
$$
Thus the law at a fixed time $t$ is uniquely defined. A standard argument allows 
to prove that this implies uniqueness for the martingale problem.  $\Box$

For $d=3$ the proof of uniqueness still works. The problem is that we cannot prove existence of a solution of the martingale problem. More precisely, we cannot prove Proposition \ref{p3.2}.

We can prove existence and uniqueness in $d=3$ for the the martingale problem where $\mathcal H$ is replaced by $\mathcal G_1$, but since the definition of $\mathcal G_1$ depends on the semigroup, this does not give any real information.

We have the following weaker result on a stopped martingale problem.
\begin{Definition}
\label{d2}
We say that a probability measure $\P_x$  on $C([0,T];D(A))$ is a solution of the stopped martingale problem
associated to \eqref{e1.1}  if
$$
\P_x(\eta(0)=1)=1,
$$
and for any $f\in \widetilde{\mathcal H}_k$
$$
f(\eta(t\wedge \tau^*))-\int_0^{t\wedge \tau^*} Lf(\eta(s))ds,
$$
is a martingale with respect to the natural  filtration and
$$
\eta(t)=\eta(\tau^*),\; t\ge \tau^*.
$$ 
The stopping time $\tau^*$ is defined
by
$$
\tau^* =Ê\lim_{R\to\infty } \tau_R,\; \tau_R= \inf\{t\in [0,T],\; |A\eta(t)|\ge R\}.
$$
\end{Definition}
\begin{Theorem}
\label{t4.2}
For any $x\in D(A)$, there exists a unique solution to the stopped martingale problem.
\end{Theorem}
{\bf Proof.} 
Existence of a solution for this martingale problem is classical. A possible proof follows the same line as 
the proofs of Proposition \ref{p3.2} and Theorem \ref{t3.2}, 
see also Remark 3.5 (see also \cite{flandoli-cetraro} for more details). In fact, we may choose the
Markov family $((X(t,x))_{t\ge 0},\Omega_x,\mathcal F_x,\P_x)$ constructed in 
\cite{DO06}. It is easy to see that $X(t,x)$ is continuous up to $\tau^*$. We slightly change
notation and set $X(t,x)=X(t\wedge\tau^*,x)$.

Uniqueness follows from a similar argument as in Theorem \ref{t4.1}. 
For $\epsilon>0$, we define $(S^\epsilon(t))_{t\ge 0}$ similarly as 
$(S_t)_{t\ge 0}$ but we replace $ e^{-K\int_{0}^{t} |A\eta(s)|^2ds} $ by 
$ e^{-\epsilon\int_{0}^{t} |A\eta(s)|^4ds} $ in \eqref{e1.4bis}. Proceeding as above, we then 
define $N_\epsilon$, $\bar N_\epsilon$, $\widetilde{\mathcal H}_k^\epsilon$, and prove 
that $\widetilde{\mathcal H}_k^\epsilon$ is a $\pi$-core for $\bar N_\epsilon$ and
$N_\epsilon \varphi = \bar N_\epsilon \varphi$ for $\varphi \in \widetilde{\mathcal H}_k^\epsilon$.

Let $\tilde\P_x$ be a solution to the martingale problem and $f\in {\mathcal E}_3^k$. For $\lambda,\; 
\epsilon >0$, 
we set $\varphi = (\lambda - \bar N_\epsilon)^{-1}$, then $\varphi \in  \widetilde{\mathcal H}_k^\epsilon$.

By It\^o formula - note that in Definition \ref{d2} it is required that the measure is supported by 
$C([0,T];D(A))$ - we prove that 
$$
\begin{array}{l}
\ds e^{-\epsilon\int_{0}^{t} |A\eta(s)|^4ds} \varphi(\eta(t) )-  \int_{0}^{t} \left( -\epsilon|A\eta(s)|^4 \varphi(\eta(s)) + L \varphi(\eta(s)) \right) e^{-\epsilon\int_{0}^{s} |A\eta(\sigma)|^4d\sigma}ds\\
\\
\ds = e^{-\epsilon\int_{0}^{t} |A\eta(s)|^4ds} \varphi(\eta(t) )-  \int_{0}^{t} N_\epsilon \varphi(\eta(s)) e^{-\epsilon\int_{0}^{s} |A\eta(\sigma)|^4d\sigma}ds
\end{array}
$$
is also a martingale. We have used :
$$
e^{-\epsilon\int_{0}^{t} |A\eta(s)|^4ds}=0,\; t\ge \tau^*.
$$
Thus:
$$
\tilde \E_x\left(e^{-\epsilon\int_{0}^{t} |A\eta(s)|^4ds} \varphi(\eta(t) )-  \int_{0}^{t} N_\epsilon \varphi(\eta(s)) e^{-\epsilon\int_{0}^{s} |A\eta(\sigma)|^4d\sigma}ds\right)= \varphi(x).
$$
We multiply by $e^{-\lambda t}$ and integrate over $[0,\infty)$ and obtain, since $\bar N_\epsilon \varphi = N_\epsilon\varphi$,
$$\begin{array}{l}
\ds\tilde\E_x\left( \int_0^\infty e^{-\lambda t -\epsilon\int_{0}^{t} |A\eta(s)|^4ds} f(\eta(t)) dt \right)\\
\\
\ds
=  \varphi(x)= (\lambda - \bar N_\epsilon)^{-1}f (x)= \int_0^\infty e^{-\lambda t}S_t^\epsilon f(x) dt. 
\end{array}
$$
By dominated convergence, we may let $\epsilon \to 0$ and obtain
$$
\tilde\E_x\left( \int_0^\infty e^{-\lambda t } 1\hspace{-0,1cm}{\rm I}_{t\le \tau^*} f(\eta(t)) dt \right)
=   \int_0^\infty e^{-\lambda t}S_t^0 f(x) dt,
$$
where $S_t^0f(x) = \lim_{\epsilon\to 0} S_t^\epsilon f(x) = \E_x( 1\hspace{-0,1cm}{\rm I}_{t\le \tau^*} f(X(t,x)))$.
The conclusion follows. $\Box$

\section{Technical results}
In all this section, we assume that $d=2$. Also, for $s\in \R$, we set $|\cdot|_s=|(-A)^s \cdot|$.
\begin{Lemma}
\label{l1}
There exists $c$ depending on $T,Q,A$ such that
$$
\E\left(\sup_{t\in[0,T]} |X(t,x)|^2+\int_0^T|X(s,x)|_1^2ds  \right)\le c(1+|x|^2),
$$
$$
\E\left(\sup_{t\in[0,T]} |X(t,x)|^4+\int_0^T|X(s,x)|^2\;|X(s,x)|_1^2ds  \right)\le c(1+|x|^4).
$$
\end{Lemma}
{\bf Proof}. We first apply It\^o's formula to $\frac12\;|x|^2$ (as usual the computation is formal and it should be justified by Galerkin approximations):
$$
\frac12\;d|X(t,x)|^2+|X(t,x)|^2_1dt=(X(t,x),\sqrt Q dW)+\frac12\;\;\mbox{\rm Tr}\;Qdt.
$$
We deduce, thanks to a classical martingale inequality,
$$
\begin{array}{l}
\ds\E\left(\frac12\;\sup_{t\in[0,T]} |X(t,x)|^2+\int_0^T|X(s,x)|_1^2ds  \right)
\\
\\
\ds\le\E\left(\sup_{t\in[0,T]}\left|\int_0^t(X(s,x),\sqrt Q dW(s))\right|  \right)+\frac12\;(|x|^2+\mbox{\rm Tr}\;Q \;T)\\
\\
\ds\le 2\E\left(   \left(\int_0^T|\sqrt Q X(s,x)|^2ds  \right)^{1/2}\right)+\frac12\;(|x|^2+\mbox{\rm Tr}\;Q \; T)\\
\\
\ds \le \frac12\;\E\int_0^T|X(s,x)|^2ds+C+\frac12\;|x|^2,
\end{array}
$$
where $C$ depends on $T,Q,A$. It follows that
\begin{equation}
\label{ea}
\E\left(\sup_{t\in[0,T]} |X(t,x)|^2+\int_0^T|X(s,x)|_1^2ds  \right)\le C+|x|^2.
\end{equation}
We now apply It\^o's formula to $\frac14\;|x|^4,$
$$
\begin{array}{l}
\ds\frac14\;d|X(t,x)|^4+|X(t,x)|^2\;|X(t,x)|^2_1dt=|X(t,x)|^2(X(t,x),\sqrt Q dW)\\
\\
\ds+\left(\frac12\;\mbox{\rm Tr}\;Q|X(t,x)|^2+|\sqrt Q X(t,x)|^2\right)dt\\
\\
\ds\le |X(t,x)|^2(X(t,x),\sqrt Q dW)+c|X(t,x)|^2dt.
\end{array}
$$
We deduce
$$
\begin{array}{l}
\ds\E\left(\frac14\;\sup_{t\in[0,T]} |X(t,x)|^2+\int_0^T|X(s,x)|^2\;|X(s,x)|_1^2ds  \right)
\\
\\
\ds\le\E\left(\sup_{t\in[0,T]}\left|\int_0^t|X(s,x)|^2(X(s,x),\sqrt Q dW(s))\right|  \right)\\
\\
\ds+c\E\int_0^T|X(s,x)|^2ds+\frac14\;|x|^4\\
\\
\ds\le 2\E\left(\left(\int_0^T|X(s,x)|^4\;|\sqrt QX(s,x)|^2ds   \right)^{1/2}     \right)  +c(1+|x|^4)
\\
\\
\ds  \le 2\E\left(\sup_{t\in[0,T]}|X(s,x)|^2
\left(\int_0^T|\sqrt QX(s,x)|^2ds   \right)^{1/2}    \right) +c(1+|x|^4)\\
\\
\ds\le \frac18\;\E\left(\sup_{t\in[0,T]} |X(t,x)|^4\right)
+c\E\left(\int_0^T|\sqrt QX(s,x)|^2 ds  \right)+c(1+|x|^4).
\end{array}
$$
Since $\sqrt Q$ is a bounded operator, using \eqref{ea} we deduce
$$
\E\left(\sup_{t\in[0,T]} |X(t,x))|^4+\int_0^T|X(s,x)|^2\;|X(s,x)|_1^2ds\right)\le(1+|x|^4).
$$
$\Box$
\begin{Lemma}
\label{l2}
There exists $c$ depending on $T,Q,A$ such that
$$
\begin{array}{l}
\E\left(\sup_{t\in[0,T]}e^{-c\int_0^t|X(s,x)|^2\;|X(s,x)|_1^2ds} |X(t,x)|^2_1\right)\\
\\
\ds+\E\left(\int_0^Te^{-c\int_0^s|X(\sigma,x)|^2\;|X(\sigma,x)|_1^2d\sigma}|X(s,x)|_2^2ds\right)\le c(1+|x|_1^2).
\end{array}
$$
\end{Lemma}
{\bf Proof}. We apply It\^o's  formula to
$$
e^{-c\int_0^t|X(s,x)|^2\;|X(s,x)|_1^2ds}|X(t,x)|^2_1,
$$
and obtain
$$
\begin{array}{l}
\ds\frac12\;d\left(e^{-c\int_0^t|X(s,x)|^2\;|X(s,x)|_1^2ds}|X(t,x)|^2_1\right)+e^{-c\int_0^t|X(s,x)|^2\;|X(s,x)|_1^2ds}\;|X(t,x)|_2^2 dt
\\
\\
\ds=e^{-c\int_0^t|X(s,x)|^2\;|X(s,x)|_1^2ds}\left(-c
|X(t,x)|^2\;|X(t,x)|_1^4+(b(X(t,x)),AX(t,x))\right)dt\\
\\
\ds+(Ax,\sqrt Q dW)-\frac12\;\mbox{\rm Tr}\;[AQ]dt.
\end{array}
$$
We have
$$
\begin{array}{l}
(b(x),Ax)\le |b(x)|\;|Ax|
\\
\\
\le\tilde c|x|_{L^4}\;|\nabla x|_{L^4}\;|Ax|\\
\\
\le\tilde c|x|^{1/2}\;|x|_{1}\;|x|_2^{3/2}\\
\\
\le\frac12\;|x|_2^2+\tilde c|x|^2\;|x|^4_{1}.
\end{array}
$$
We deduce that if $c\ge \tilde c$,
$$
\begin{array}{l}
\ds\frac12\;d\left(e^{-c\int_0^t|X(s,x)|^2\;|X(s,x)|_1^2ds}|X(t,x)|^2_1\right)+\frac12\;e^{-c\int_0^t|X(s,x)|^2\;|X(s,x)|_1^2ds}\;|X(t,x)|_2^2 dt
\\
\\
\ds\le e^{-c\int_0^t|X(s,x)|^2\;|X(s,x)|_1^2ds} (AX(t,x),\sqrt Q dW)+cdt
\end{array}
$$
and
$$
\begin{array}{l}
\ds\E\left(\sup_{t\in[0,T]}e^{-c\int_0^t|X(s,x)|^2\;|X(s,x)|_1^2ds}|X(t,x)|^2_1\right)\\
\\
\ds+\E\left(\int_0^Te^{-c\int_0^s|X(\sigma,x)|^2\;|X(\sigma,x)|_1^2d\sigma}\;|AX(s,x)|^2 ds\right)
\\
\\
\ds\le 2\E\left(\left(\int_0^Te^{-2c\int_0^s|X(\sigma,x)|^2\;|X(\sigma,x)|_1^2d\sigma}\;|\sqrt QX(s,x)|^2 ds\right)^{1/2}\right)\\
\\
+cT+|x|^2_1.
\end{array}
$$
Since Tr $(QA)<\infty$, we know that $QA$ is a bounded operator and
$$
\begin{array}{l}
\ds\E\left(\left(\int_0^Te^{-2c\int_0^s|X(\sigma,x)|^2\;|X(\sigma,x)|_1^2d\sigma}\;|\sqrt QX(s,x)|^2 ds\right)^{1/2}\right)\\
\\
\ds\le c\E\left(\left(\int_0^T|X(s,x)|_1^2 ds\right)^{1/2}\right)\\
\\
\le (|x|+1),
\end{array}
$$
by Lemma \ref{l1}. The result follows. $\Box$
\begin{Lemma}
\label{l2bis}
For any $k\in \N$, there exists $c$ depending on $k,T,Q,A$ such that
$$
\E\left(\sup_{t\in[0,T]}e^{-c\int_0^t|X(s,x)|^2\;|X(s,x)|_1^2ds} |X(t,x)|^k_1\right)
\le c(1+|x|_1^k).
$$
\end{Lemma}
The proof of this Lemma follows the same argument as above. It is left to the reader.

\begin{Proposition}
\label{p3.1}
For any $k\in \N$, $\epsilon >0$, there exists $C(\epsilon,k,T,Q,A)$ such that for any $m\in\N$,
$x\in D(A)$, $t\in [0,T]$,
$$
\E\left(e^{-\epsilon\int_0^t|(-A)^{1/2}X_m(s,x)|^6 ds}|AX_m(t,x)|^k\right) \le C(\epsilon,k,T,Q,A)(1+|Ax|^k).  
$$
\end{Proposition}
{\bf Proof.}
Let us set
$$
Z(t)=\int_0^t e^{(t-s)A}\sqrt QdW(s),\quad Y(t)=X(t,x)-z(t).
$$
Then, by the factorization method (see \cite{DPZ}),
\begin{equation}
\label{eabis}
\E\left(\sup_{t\in[0,T]} |Z(t)|^h_{2+\epsilon}  \right)\le C, 
\end{equation}
for any $\epsilon<g$, and we have
$$
\frac{dY}{dt}=AY+b(Y+Z).
$$
We take the scalar product with $A^2Y$:
$$
\frac12\;\frac{d}{dt}|Y|_2^2+|Y|_3^2=(b(Y+Z),A^2Y)=((-A)^{1/2}b(Y+Z),A^{3/2}Y).
$$ 
We have
$$
\begin{array}{l}
\ds|(-A)^{1/2}b(Y+Z)|=|\nabla b(Y+Z)|
\\
\\
\le c\left( |Y+Z|^2_{W^{1,4}} +|Y+Z|_{L^p} |Y+Z|_{W^{2,q}}\right) ,
\end{array}
$$
where $\frac1p+\frac1q=\frac12.$

By Gagliardo-Nirenberg inequality
$$
|Y+Z|^2_{W^{1,4}}\le c|Y+Z|_1\;|Y+Z|_2.
$$
Setting $\frac1p=\frac12-\frac{s}2$ we have by Sobolev's embedding
$$
|Y+Z|_{L^p}\;|Y+Z|_{W^{2,q}}\le c|Y+Z|_{s/2}\;|Y+Z|_{3-s/2}.
$$
Therefore
$$
\begin{array}{l}
\ds ((-A)^{1/2}b(Y+Z),((-A)^{3/2}Y)\le c|Y+Z|_{1}\;|Y+Z|_{2}\;|Y|_{3}
\\
\\
\le c|Y+Z|_{s/2}\;|Y+Z|_{3-s/2}\;|Y|_{3}
\\
\\
\ds\le \frac14\;|Y|_{3}^2+c|Y+Z|_{1}^2\;|Z|_2^2+c|Y+Z|_{1}^2\;|Y|_2^2\\
\\
+c |Y+Z|_{s/2}\;|Y|_{3-s/2}\;|Y|_{3}+c|Y+Z|^2_{s/2}\;|Z|^2_{3-s/2}.
\end{array}
$$
Since
$$
\begin{array}{l}
\ds |Y+Z|_{s/2}\;|Y|_{3-s/2}\;|Y|_{3}\le
|Y+Z|_{s/2}\;|Y|_1^{s/4}\;|Y|_{3}^{2-s/4}
\\
\\
 \le c|Y+Z|^{8/s}_{s/2}\;|Y|_1^{2}+\frac14\;|Y|_{3}^2,
\end{array}
$$
we finally get
$$
\begin{array}{l}
\ds \frac{d}{dt}\;|Y|_2^2\le c|Y+Z|^2_{1}\;|Y|^2_{2}
\\
\\
+c\left(|Y+Z|^2_{1}\;|Z|^2_{2}+|Y+Z|^2_{s/2}\;|Z|^2_{3-s/2}+|Y+Z|^{8/s}_{s/2}\;|Y|^{2}_{1}   \right) 
\end{array}
$$
and
$$
\begin{array}{l}
\ds |Y(t)|^2_2\le e^{c\int_0^t|Y+Z|^2_{1}ds}
\\
\\
\ds \left(|x|_2^2+c\int_0^t(|Y+Z|^2_{1}\;|Z|_2^2+|Y+Z|^2_{s/2}\;|Z|^2_{3-s/2}+|Y+Z|^{8/s}_{1}\;|Y|^{2}_{1/2})ds\right).
\end{array}
$$
We then write by H\"older and Poincar\'e inequalities
$$
\begin{array}{l}
\ds e^{-\epsilon\int_0^t|Y+Z|^6_{1}ds}|Y(t)|_2^k\le c_k
e^{-\epsilon\int_0^t|Y+Z|^6_{1}ds+c_k\int_0^t|Y+Z|^2_{1}ds}
\\
\\
\times\left(|x|_2^k+ \int_0^t(|Y+Z|^4_{1}+|Z|^4_{2+g} +|Y+Z|^{16/s}_{1})ds \right)
\end{array}
$$
(we choose $3-s/2<2+g$ and set $\epsilon=1-s/2$).

The conclusion follows from Lemma \ref{l2bis} and by the boundedness of 
$x\mapsto -\epsilon x^6+cx^4$. $\Box$


\newpage
\footnotesize
 \begin{thebibliography}{9}


\bibitem{kupiainen} J. Bricmont, A. Kupiainen
and R. Lefevere, {\it
Exponential mixing for the $2D$ Navier-Stokes dynamics}, Comm. Math. Phys., {\bf 230}, no 1, 87-132,
2002.


\bibitem{chueschov-kuksin} I. Chueshov, S. B. Kuksin, {\it On the random kick-forced 
3D Navier-Stokes equations in a thin domain},  Arch. Ration. Mech. Anal. 188, no. 1, 117--153, 2008.

\bibitem{chueschov-kuksin2} I. Chueshov, S. B. Kuksin, {\it Stochastic 3D Navier-Stokes equations in a thin domain and its $\alpha$-approximation}, Preprint (2007).
 
\bibitem{DPD-NS3D} G. Da Prato,  A. Debussche,  {\it Ergodicity for the
3D stochastic Navier-Stokes
equations,}  J. Math. Pures et Appl., 82,  877-947, 2003.

\bibitem{DPZ} G. Da Prato,  J. Zabczyk, {\it Stochastic equations
in infinite
dimensions,}
   Encyclopedia of Mathematics and its Applications,  Cambridge University
Press, 1992.


\bibitem{DPZ2} G. Da Prato and J. Zabczyk,   {\it Ergodicity for
infinite dimensional
systems},    London  Mathematical Society Lecture Notes, {\bf 229}, Cambridge University
Press, 1996.

\bibitem{DO06} A. Debussche and C. Odasso, {\it  Markov solutions for  the $3D$ stochastic
Navier--Stokes equations with state dependent noise}, Journal of Evol. Equ., vol 6, no 2, 305-324,
2006.


\bibitem{mattingly} W. E, J.C. Mattingly, Y. G. Sinai,
{\it Gibbsian dynamics and ergodicity for the stochastically forced Navier-Stokes equation},
Commun. Math. Phys. {\bf 224}, 83--106, 2001.


\bibitem{flandoli-cetraro}  F. Flandoli,  {\it An introduction to 3D stochastic Fluid Dynamics}, CIME
lecture notes 2005.



\bibitem{flandoli-maslowski} F. Flandoli and B. Maslowski, {\it Ergodicity of
the  2D Navier-Stokes equations under random perturbations}, Comm. Math. Phys., {\bf 172},
no 1, 119-141, 1995.


\bibitem{romito-flandoli-selection} F. Flandoli, M. Romito, {\it Markov selections for 
the 3D stochastic Navier-Stokes equation}, Probab. Theory Relat. Fields 140, nr. 3-4 (2008), 407--458.

\bibitem{GRZ} B. Goldys, M. R\"ockner and X. Zhang, {\it Martingale solutions and Markov selections for stochastic partial differential equations}, Bielefeld, Preprint, 2008.


\bibitem{HM} M. Hairer, J.C. Mattingly, {\it Ergodicity of the 2D Navier-Stokes
Equations with Degenerate Stochastic Forcing},  Annals of Mathematics, vol. 164 no. 3,   2006. 

\bibitem{HM2} M. Hairer, J.C. Mattingly. {\it   Spectral gaps in Wasserstein distances and the 2D stochastic Navier-Stokes equations}. 2006. Preprint. 

\bibitem{kuksinBook} S B Kuksin, {\it Randomly forced nonlinear PDEs and statistical hydrodynamics in 
2 space dimensions}, Europear Mathematical Society Publishing House, 2006. 

\bibitem{kuksin} S. Kuksin and  A. Shirikyan, {\it Ergodicity for the radomly forced
2D Navier-Stokes equations}, Math. Phys. Anal. and Geom., {\bf 4}, no 2, 147-195, 2001.

\bibitem{kuksin2} S. Kuksin and  A. Shirikyan, {\it A coupling approach to randomly forced 
nonlinear PDEs. I.}, Comm. Math. Phys., {\bf 221}, 351-366, 2001.

\bibitem{mattingly2} J.C. Mattingly, {\it Exponential convergence for the stochastically forced Navier-Stokes
equations and other partially dissipative dynamics}. Commun. Math. Phys. {\bf 230},  no. 3, 421-462, 2002.
 
\bibitem{mattingly3}  J. C. Mattingly and E. Pardoux {\it Malliavin calculus for the stochastic 
2D Navier-Stokes equation},  Comm. Pure Appl. Math, {\bf 59},  no. 12   1742-1790,
2006. 

\bibitem{priola}  E. Priola, {\it On a class of Markov type semigroups in spaces of uniformly continuous
 and bounded functions}, Studia Math. 136 (1999), no. 3, 271--295.
 
\bibitem{roeckner-sobol} M. R\"ockner, Z. Sobol, {\it A new approach to Kolmogorov equations in infinite dimensions and applications to the stochastic 2d Navier Sokes equation}  C. R. Acad. Sci. Paris, S\'erie I, 4 pp., 2007. 
 
 \bibitem{romito} M. Romito, {\it Analysis of equilibrium states of Markov solutions to the 3D Navier-Stokes equations driven by additive noise}, to appear on J. Stat. Phys.

\bibitem{Odasso} C. Odasso,
{\it Exponential mixing for the $3$D Navier-Stokes equation}. 2005. To appear in Commun. Math. Phys.

\bibitem{shirikyan3D} A. Shirikyan, {\it Qualitative properties of stationary measures for three-dimensional Navier-Stokes equations}, J. Funct. Anal. 249 (2007), 284-306.


\end{thebibliography}
\end {document}